\documentclass[12pt]{amsart}

\usepackage{latexsym}
\usepackage{amssymb}

\newtheorem{thm}{Theorem}[section]
\newtheorem{lemma}[thm]{Lemma}

\newtheorem{cor}[thm]{Corollary}
\newtheorem{prop}[thm]{Proposition}

\newtheorem{Example}[thm]{Example}
\newtheorem{Remark}[thm]{Remark}
\newtheorem{Alg}[thm]{Algorithm}
\newtheorem{Defn}[thm]{Definition}

\newenvironment{remark}{\begin{Remark}\rm}
		{\mbox{}~\hfill$\square$\end{Remark}}

\newenvironment{defn}{\begin{Defn}\rm}{\end{Defn}}

	{\begin{list}
		{\noindent\makebox[0mm][r]{(\roman{enumi})}}
		{\usecounter{enumi} \topsep=1.5mm \itemsep=0mm}
		
	}
	{\end{list}}

	{\begin{list}
		{\noindent\makebox[0mm][r]{\arabic{enumi}.}}
		{\usecounter{enumi} \topsep=1.5mm \itemsep=-1mm}
	}
	{\end{list}}

\def\bem#1{\textbf{#1}}

\def\<{\langle}
\def\>{\rangle}
\def\0{{\mathbf 0}}
\def\1{{\mathbf 1}}
\def\NN{{\mathbb N}}
\def\RR{{\mathbb R}}
\def\WW{{\hspace{.1ex}\overline{\hspace{-.1ex}W\hspace{-.1ex}}\hspace{.1ex}}{}}
\def\ZZ{{\mathbb Z}}
\def\mm{{\mathfrak m}}
\def\xx{{\mathbf x}}

\def\th{{\rm th}}

\def\hom{{\rm Hom}}
\def\tot{{\rm tot\hspace{.2ex}}}
\def\char{{\rm char}}
\def\link{{\rm link}}
\def\rays{{\rm rays}}
\def\facets{{\rm facets}}

\def\dv{{\Delta^{\!\star}}}

\def\eext{{\hspace{.6pt}{\underline{\hspace{-.6pt}{\rm Ext}\hspace{-.6pt}}%
				\hspace{.6pt}}{\scriptscriptstyle \,}}}

\def\hhom{{\hspace{.6pt}{\underline{\hspace{-.6pt}{\rm Hom}\hspace{-.6pt}}%
				\hspace{.6pt}}{\scriptscriptstyle \,}}}
\def\into{\hookrightarrow}
\def\spot{{\hbox{\raisebox{1.7pt}{\large\bf .}}\hspace{-.5pt}}}

\def\minus{\smallsetminus}
\def\congto{\stackrel{\begin{array}{@{}c@{\;}}
                        \\[-4ex]\scriptstyle\approx\\[-1.6ex]
                      \end{array}}\to}
\def\implies{\Rightarrow}
\def\nothing{\varnothing}
\def\fillrightmap{\mathord- \mkern-6mu
	\cleaders\hbox{$\mkern-2mu \mathord- \mkern-2mu$}\hfill
	\mkern-6mu \mathord\rightarrow}

\def\wt#1{{\widetilde{#1}}}
\def\mod#1{\ ({\rm mod}\:#1)}
\def\text#1{\hbox{\scriptsize #1}}
\newcommand{\aoverb}[2]{{\genfrac{}{}{0pt}{1}{#1}{#2}}}
\def\twoline#1#2{{\scriptstyle \aoverb{{#1}_{\,}}{{#2}^{\,}}}}

\begin{document}

\flushbottom

\title[Cohen--Macaulay quotients via irreducible resolutions]
	{Cohen--Macaulay quotients of normal semigroup rings via
	irreducible resolutions}
\renewcommand{\thefootnote}{\fnsymbol{footnote}}
\author[Ezra Miller]{Ezra Miller$^*$}%
%

\begin{abstract}
\noindent
For a radical monomial ideal $I$\/ in a normal semigroup ring~$k[Q]$,
there is a unique minimal \bem{irreducible resolution} $0 \to k[Q]/I \to
\WW^0 \to \WW^1 \to \cdots$\ by modules $\WW^i$ of the form $\bigoplus_j
k[F_{ij}]$, where the $F_{ij}$ are (not necessarily distinct) faces
of~$Q$.  That is, $\WW^i$ is a direct sum of quotients of $k[Q]$ by prime
ideals.  This paper characterizes Cohen--Macaulay quotients $k[Q]/I$ as
those whose minimal irreducible resolutions are \bem{linear}, meaning
that $\WW^i$ is pure of dimension $\dim(k[Q]/I) - i$ for $i \geq 0$.  The
proof exploits a graded ring-theoretic analogue of the Zeeman spectral
sequence~\cite{zeeIII}, thereby also providing a combinatorial
topological version involving no commutative algebra.  The
characterization via linear irreducible resolutions reduces to the
Eagon--Reiner theorem \cite{ER} by Alexander duality
\mbox{when $Q = \NN^d$.}
\vskip 1ex
\noindent
{{\it 2000 AMS Classification:} 13C14, 14M05, 13D02, 55Txx (primary)
14M25, 13F55 (secondary)}
\end{abstract}

\maketitle

\footnotetext{$^*$The author was funded by the National Science
Foundation.}
\renewcommand{\thefootnote}{\arabic{footnote}}

{}


\section{Introduction}

\label{sec:intro}

Let $Q \subseteq \ZZ^d$ be a normal affine semigroup.  We assume for
simplicity that~$Q$ generates~$\ZZ^d$ as a group, and that~$Q$ has
trivial unit group.  The real cone $\RR_{\geq 0} Q$
is a polyhedral cell complex.  Endow it with an incidence function
$\varepsilon$, and let $\Delta \subseteq \RR_{\geq 0} Q$ be a closed
polyhedral subcomplex.  Corresponding to $\Delta$ is the ideal $I_\Delta$
inside the semigroup ring~$k[Q]$, generated (as a $k$-vector space) by
all monomials in $k[Q]$ not lying on any face of~$\Delta$.  Thus
$k[Q]/I_\Delta$ is spanned by monomials lying in~$\Delta$.

This paper has three goals:
\begin{itemize}
\item
Define the notion of \bem{irreducible resolution} for $Q$-graded
$k[Q]$-modules.

\item
Introduce the \bem{Zeeman double complex} for $\Delta$.

\item
Characterize Cohen--Macaulay quotients $k[Q]/I_\Delta$ in terms of the
above items.
\end{itemize}

An irreducible resolution (Definition~\ref{d:irr}) of a $\ZZ^d$-graded
$k[Q]$-module is an injective-like resolution, in which the summands are
quotients of $k[Q]$ by irreducible monomial ideals rather than
indecomposable injectives.  Minimal irreducible resolutions exist
uniquely up to isomorphism for all $Q$-graded modules~$M$
(Theorem~\ref{t:irr}).  When $M = k[Q]/I_\Delta$, every summand is
isomorphic to a semigroup ring~$k[F]$, considered as a quotient module
of~$k[Q]$, for some face $F \in \Delta$ (Corollary~\ref{c:irr}).

The Zeeman double complex $D(\Delta)$
consists of $k[Q]$-modules that are direct sums of semigroup rings $k[F]$
for faces $F \in \Delta$ (Definition~\ref{d:zeeman}).  Its naturally
defined differentials come
from the incidence function on~$\Delta$.

Here is the idea behind the Cohen--Macaulay criterion,
Theorem~\ref{t:cm}.  Although the total complex of the Zeeman double
complex $D(\Delta)$ is an example of an irreducible resolution
(Proposition~\ref{p:tot}), its large number of summands keeps it far from
being minimal.  However, the cancellation afforded by the horizontal
differential of $D(\Delta)$ sometimes causes the resulting vertical
differential (on the horizontal cohomology) to be a minimal irreducible
resolution.  This fortuitous cancellation occurs precisely when $\Delta$
is Cohen--Macaulay over~$k$, in which case the horizontal cohomology
occurs in exactly one column of $D(\Delta)$.

Part~\ref{ordinary} of Theorem~\ref{t:cm}, which characterizes the
Cohen--Macaulay property by collapsing at $E^1$ of the ordinary Zeeman
spectral sequence for $\Delta$ (Definition~\ref{d:sequence}), may be of
interest to algebraic or combinatorial topologists.  Its statement as
well as its proof
are independent of the surrounding commutative algebra.

The methods involving Zeeman double complexes and irreducible resolutions
should have applications beyond those investigated here; see
Section~\ref{sec:further} for possibilities.

\subsection*{Notational conventions}
By the assumptions on the semigroup~$Q$ in the first paragraph above, $Q$
is the intersection with $\ZZ^d$ of the positive half-spaces defined by
primitive integer-valued functionals $\tau_1, \ldots, \tau_n$ on~$\ZZ^d$.
In particular the $i^\th$ \bem{facet} of~$Q$ (for $i = 1,\ldots,n$) is
the subset $F_i \subseteq Q$ on which $\tau_i$ vanishes.  More generally,
an arbitrary \bem{face} of~$Q$ is defined by the vanishing of a linear
functional on $\ZZ^d$ that is nonnegative on~$Q$.  The (Laurent) monomial
in~$k[\ZZ^d]$ with exponent~$\alpha$ is denoted by $\xx^\alpha$, although
sets of monomials in $k[\ZZ^d]$ are frequently identified with their
exponent sets in~$\ZZ^d$.

All cellular homology and cohomology groups are taken with coefficients
in the field~$k$ unless otherwise stated.
We work here always with nonreduced (co)homology of the usually unbounded
polyhedral complex~$\Delta$, which corresponds to the reduced
(co)homology of an always bounded transverse hyperplane section
of~$\Delta$, homologically shifted by~$1$.

All modules in this paper, including injective modules, are
$\ZZ^d$-graded unless otherwise stated.  Elementary facts regarding the
category of $\ZZ^d$-graded $k[Q]$-modules, especially $\ZZ^d$-graded
injective modules, hulls, and resolutions,
can be found in~\cite{GWii}.

\section{Irreducible resolutions}

\label{sec:irr}

Recall that an ideal $W$ inside of $k[Q]$ is \bem{irreducible} if
$W$ can't be expressed as an intersection of two ideals properly
containing~it.

\begin{defn} \label{d:irr}
An \bem{irreducible resolution} $\WW^\spot$ of a $k[Q]$-module~$M$ is an
exact sequence
$$
  0 \to M \to \WW^0 \to \WW^1 \to \cdots \qquad\quad \WW^i = \bigoplus_{j
  = 1}^{\mu_i} k[Q]/W^{ij}
$$
in which each $W^{ij}$ is an irreducible ideal of~$k[Q]$.  The
irreducible resolution is called \bem{minimal} if all the numbers $\mu_i$
are simultaneously minimized (among irreducible resolutions of~$M$), and
\bem{linear} if $\WW^i$ is pure of Krull dimension $\dim(M) - i$ for
all~$i$.  (By convention, modules of negative dimension are zero.)%
\end{defn}

The fundamental properties of quotients $\WW := k[Q]/W$ by irreducible
monomial ideals $W$ are inherited from the corresponding properties of
indecomposable injective modules.  Recall that each such
\bem{indecomposable injective module} is a vector space $k\{\alpha +
E_F\}$ spanned by the monomials in $\alpha + E_F$, where
\begin{eqnarray} \label{eq:EF}
  E_F &=& \{f-a \mid f \in F \hbox{ and } a \in Q\}
\end{eqnarray}
is the negative tangent cone along the face $F$ of~$Q$.  The vector space
$k\{\alpha+E_F\}$ carries an obvious structure of $k[Q]$-module.  In what
follows, the $\ZZ^d$-graded injective hull of a $\ZZ^d$-graded module~$M$
\cite{GWii} is denoted by~$E(M)$, so that, in particular, $E(k[F]) =
k\{E_F\}$.  Define the \bem{$Q$-graded part} of~$M$ to be the submodule
$\bigoplus_{a \in Q} M_a$ generated by elements whose degrees lie in~$Q$.

\begin{lemma} \label{l:irrideal}
A monomial ideal $W$ is irreducible
if and only if $\WW := k[Q]/W$ is the $Q$-graded part of some
indecomposable injective module.
\end{lemma}
\begin{proof}
($\Leftarrow$) The module $k\{\alpha + E_F\}_Q$ is clearly isomorphic
to~$\WW$ for some ideal~$W$.  Supposing that $W \neq k[Q]$, we may as
well assume $\alpha \in Q$ by adding an element far inside~$F$, so that
$\xx^\alpha \in \WW$ generates an essential submodule $k\{\alpha + F\}$.
Suppose $W = W_1 \cap W_2$.  The copy of $k\{\alpha+F\}$ inside $\WW$
must include into $\WW_{\!j}$ for $j = 1$ or $2$.  Indeed, if both
induced maps $k\{\alpha+F\} \to \WW_{\!j}$ have nonzero kernels, then
they intersect in a nonzero submodule of $k\{\alpha+F\}$ because $k[F]$
is a domain.  The essentiality of $k\{\alpha+F\} \subseteq \WW$ then
forces $\WW \to \WW_{\!j}$ to be an inclusion for some~$j$.  We conclude
that $W$ contains this $W_j$, so $W = W_j$ is irreducible.

($\implies$) Let $W$ be an irreducible ideal and $\WW = k[Q]/W$.
Cconsidering the injective hull~$E(\WW) = J_1 \oplus \cdots \oplus J_r$,
the composite $k[Q] \to \WW \to E(\WW)$ has kernel $W = W_1 \cap \cdots
\cap W_r$, where $\WW_{\!j} = (J_j)_Q$.  Since $W$ is irreducible, we
must have $W_j = W$ for some~$j$.  We conclude that $E(\WW) = J_j$, and
$W = W_j$.
\end{proof}

\begin{lemma} \label{l:isom}
For any finitely generated
module~$M$, there exists $\beta \in \ZZ^d$ such that $M_\beta \neq 0$,
and for all $\gamma \in \beta + Q$, the inculsion $M_\gamma \into
E(M)_\gamma$ is an isomorphism.
\end{lemma}
\begin{proof}
Suppose that $E(M) = \bigoplus_{\alpha,F} k\{\alpha +
E_F\}^{\mu(\alpha,F)}$, where we assume that $\alpha + E_F \neq \alpha' +
E_{F'}$ whenever $(\alpha,F) \neq (\alpha',F')$.  Now fix a pair
$(\alpha,F)$ so that $\alpha + E_F$ is maximal inside~$\ZZ^d$ among all
such subsets appearing in the direct sum.  Clearly we may assume $F$ is
maximal among faces of~$Q$ appearing in the direct sum.  Pick an element
$f$\/ that lies in the relative interior of~$F$.

By~(\ref{eq:EF}) and the maximality of~$F$, some choice of $r \in \NN$
pushes the $\ZZ^d$-degree $\alpha + r\cdot f \in \alpha + E_F$ outside of
$\alpha' + E_{F'}$ for all $(\alpha',F')$ satisfying $F' \neq F$.
Moreover, the maximality of $(\alpha,F)$
implies that $\alpha + r\cdot f \not\in \alpha' + E_F$ whenever $\alpha'
\neq \alpha$.

The prime ideal $P_F$ satisfying $k[Q]/P_F = k[F]$ is minimal over the
annihilator of~$M$.  Therefore, if $M' = (0:_M P_F)$ is the submodule
of~$M$ annihilated by~$P_F$, the composite injection $M' \into M \into
E(M)$ becomes an isomorphism onto its image after homogeneous
localization at~$P_F$---that is, after inverting the monomial $\xx^f$.
It follows that choosing $r \in \NN$ large enough forces isomorphisms
\begin{equation} \label{eq:isom}
  M_{\alpha + r\cdot f}\ \: \congto\ \: E(M)_{\alpha + r\cdot f} \ \:
  \congto\ \: \Bigl(k\{\alpha+E_F\}^{\mu(\alpha,F)}\Bigr){}_{\alpha +
  r\cdot f}\ \: \cong \ \: k^{\mu(\alpha,F)}.
\end{equation}
Setting $\beta = \alpha + r\cdot f$, the multiplication map
$\xx^{\gamma-\beta} : E(M)_\beta \to E(M)_\gamma$ for $\gamma \in \beta +
Q$ is either zero or an isomorphism, because $E(M)$ agrees with
$k\{\alpha+E_F\}^{\mu(\alpha,F)}$ in degrees $\beta$ and $\gamma$ by
construction.  The previous sentence holds with $M$ in place of $E(M)$
by~(\ref{eq:isom}), because $M$ is a submodule of $E(M)$.
\end{proof}

Not every module has an irreducible resolution, because being $Q$-graded
is a prerequisite.  However, $Q$-gradedness is the only restriction.

\begin{thm} \label{t:irr}
Let $M = M_Q$ be a finitely generated $Q$-graded module.  Then:
\begin{enumerate}
\item \label{unique}
$M$ has a minimal irreducible resolution, unique up to noncanonical
isomorphism.

\item \label{split}
Any irreducible resolution of~$M$ is (noncanonically) the direct sum of a
minimal irreducible resolution and a split exact irreducible resolution
of\/ $0$.

\item \label{finite}
The minimal irreducible resolution of $M$ has finitely many irreducible
summands in each cohomological degree.

\item \label{length}
The minimal irreducible resolution of $M$ has finite length; that is, it
vanishes in all sufficiently high cohomological degrees.

\item \label{Q-graded}
The $Q$-graded part of any injective resolution of~$M$ is an irreducible
resolution.

\item \label{injres}
Every irreducible resolution of~$M$ is the $Q$-graded part of an
injective resolution.
\end{enumerate}
\end{thm}
\begin{proof}
Lemma~\ref{l:irrideal} implies part~\ref{Q-graded}.  The remaining parts
therefore follow from part~\ref{injres} and the corresponding (standard)
facts about $\ZZ^d$-graded injective resolutions \cite{GWii}, with the
exception of part~\ref{finite}, which is false for injective resolutions
whenever $k[Q]$ isn't isomorphic to a polynomial ring.

Focusing now on part~\ref{injres}, let $\WW^\spot$ be an irreducible
resolution of~$M$, and set $J^0 = E(\WW^0)$.
The inclusion $M \into J^0$ has $Q$-graded part $M \into \WW^0$ by
Lemma~\ref{l:irrideal}.  Making use of the defining property of injective
modules, extend the composite inclusion $\WW^0/M \into \WW^1 \into
E(\WW^1)$ to a map $J^0/M \to E(\WW^1)$, and let $K^0$ be the kernel.
Then $K^0$ has zero $Q$-graded part because $\WW^0/M \into \WW^1$ is a
monomorphism.  The injective hull $K^0 \into E(K^0)$ therefore has zero
$Q$-graded part.  Extending $K^0 \into E(K^0)$ to a map $J^0/M \to
E(K^0)$ yields an injection $J^0/M \to J^1 := E(K^0) \oplus E(\WW^1)$
whose $Q$-graded part is $\WW^0/M \into \WW^1$.  Replacing~$M$,
$0$~and~$1$ by ${\rm image}(\WW^{i-1} \to \WW^i)$, $i$~and~$i+1$ in this
discussion produces the desired injective resolution by induction.%

Finally, for the length-finiteness in part~\ref{finite}, consider the set
$V(M)$ of degrees $a \in Q$ such that $M_b$ vanishes for all $b \in a+Q$.
The vector space $k\{V(M)\}$ is naturally an ideal in~$k[Q]$.
Lemma~\ref{l:isom} implies that $V(M) \subsetneq V(\WW/M)$ whenever $\WW$
is the $Q$-graded part of an injective hull of~$M$ and $M \neq 0$ (that
is, $V(M) \neq Q$).  The noetherianity of $k[Q]$ plus this strict
containment force the sequence of ideals
$$
  k\{V(M)\}\ \subseteq\ k\{V(\WW^0/M)\}\ \subseteq\ k\{V(\WW^1/{\rm
  image}(\WW^0))\}\ \subseteq\ \cdots
$$
to stabilize at the unit ideal of~$k[Q]$ after finitely many steps.%
\end{proof}

\begin{remark}
The results in this section hold just as well for unsaturated semigroups,
with the same proofs, verbatim.
\end{remark}

Examples of irreducible resolutions include Proposition~\ref{p:tot},
below, as well as the proof of Lemma~\ref{l:vert}, which contains the
irreducible resolution of the canonical module of $k[F]$ in~(\ref{eq:F}).
In general, any example of an injective resolution of any $\ZZ^d$-graded
module yields an irreducible resolution of its $Q$-graded part, although
the indecomposable injective summands with zero $Q$-graded part get
erased.  In particular, the ``cellular injective resolutions'' of
\cite{Mil2} become what should be called ``cellular irreducible
resolutions'' here.

\section{Zeeman double complex}

\label{sec:zeeman}

This section introduces the Zeeman double complex and its resulting
spectral sequences.  The total complex of the Zeeman double complex in
Proposition~\ref{p:tot} provides a natural but generally nonminimal
irreducible resolution for $k[Q]/I_\Delta$.

For each face $G \in \Delta$, let $k[G]$ be the affine semigroup ring
for~$G$, considered as a quotient of~$k[Q]$, and denote by~$e_G$ the
canonical generator of $k[G]$ in $\ZZ^d$-graded degree~$\0$.  Also, for
each face $F \in \Delta$, let $k F$ be a $1$-dimensional $k$-vector space
spanned by~$F$ in $\ZZ^d$-graded degree~$\0$.
\begin{defn} \label{d:zeeman}
Consider the $k[Q]$-module $D(\Delta) = \bigoplus_{F \supseteq G} kF
\otimes_k k[G]$ generated by
$$
  \{F \otimes_k e_G \mid F,G \in \Delta \hbox{ and }F \supseteq G\}.
$$
Doubly index the generators so that $D(\Delta)_{pq}$ is generated by
$$
  \{F \otimes e_G \mid p = \dim F \hbox{ and } -q = \dim G\},
$$
and hence $\{\0\} \otimes e_{\{\0\}} \in D(\Delta)_{00}$, with the rest
of the double complex in the fourth quadrant.  Now define the \bem{Zeeman
double complex} of $\ZZ^d$-graded $k[Q]$-modules to be $D(\Delta)$, with
vertical differential $\partial$ and horizontal differential $\delta$ as
in the diagram:
$$
\partial e_G\ =\hspace{-1ex} \sum_\twoline{G' \in G}{\text{is a facet}}
	\hspace{-1ex}\varepsilon(G',G) e_{G'}
\qquad\!
\begin{array}{@{}ccc@{}}
	F \otimes \makebox[0pt][l]{$\partial e_G$}\phantom{e_G}
\\[3pt]	\partial\uparrow\phantom{F}
\\[-5pt]F\otimes e_G & \stackrel{\textstyle \delta}{\fillrightmap}
	  & \delta F \otimes e_G
\end{array}
\quad\ 
(-1)^q \delta F\ =\hspace{-1ex} \sum_\twoline{F \in F'}{\text{is a
	facet}} \hspace{-1ex}\varepsilon(F,F') F',
$$
where the signs $\varepsilon(G',G)$ and $\varepsilon(F,F')$ come from the
incidence function on~$\Delta$.
\end{defn}

For each fixed $G$, the elements $F \otimes e_G$ generate a summand of
$D(\Delta)$ closed under the horizontal differential $\delta$.  Taking
the sum over all~$G$ yields the horizontal complex
$$
  (D(\Delta), \delta) = \bigoplus_{G \in \Delta} C^\spot(\Delta_G)
  \otimes k[G], \quad {\rm where} \quad \Delta_G = \{F \in \Delta \mid F
  \supseteq G\}
$$
is the \bem{part of $\Delta$ above}~$G$.  It is straightforward to verify
that $C^\spot(\Delta_G)$ is isomorphic to the reduced cochain complex
$\wt C^\spot \link(G,\Delta)$ of the link of~$G$ in~$\Delta$ (also known
as the vertex figure of~$G$ in~$\Delta$), but with $\nothing$ in
homological degree~$\dim G$ instead of~$-1$.

The cohomology $H^i(C^\spot(\Delta_G))$ is also called the \bem{local
cohomology~$H^i_G(\Delta)$} of $\Delta$ near~$G$.  Since the complex
$C^\spot(\Delta_G)$ is naturally a subcomplex of $C^\spot(\Delta_{G'})$
whenever $G' \subseteq G$, the natural restriction maps $H^i_G(\Delta)
\to H^i_{G'}(\Delta)$ make local cohomology into a sheaf on $\Delta$.
The following is immediate from the above discussion.

\begin{lemma} \label{l:hor}
In column~$p$, the vertical complex $(H_\delta D(\Delta), \partial)$ of
$k[Q]$-modules has $\bigoplus_{\dim G = q} k[G] \otimes_k H^p_G(\Delta)$
in cohomological degree~$-q$.  The vertical differential $\partial$ is
comprised of the natural maps $H^p_G(\Delta) \otimes e_G \to
H^p_{G'}(\Delta) \otimes \varepsilon(G',G) e_{G'}$ for \mbox{facets $G'$
of~$G$.}
\end{lemma}

We'll need to know the vertical cohomology $H_\partial D(\Delta)$ of
$D(\Delta)$, too.

\begin{lemma} \label{l:vert}
$H_\partial D(\Delta) = \bigoplus_{F \in \Delta} \omega_{k[F]}$, where
$\omega_{k[F]}$ is the canonical module of $k[F]$, and each summand
$\omega_{k[F]}$ sits along the diagonal in bidegree $(p,q) = (\dim F,
-\dim F)$.
\end{lemma}
\begin{proof}
Collecting the terms with fixed~$F$ yields the tensor product of $kF$
with
\begin{equation} \label{eq:F}
  0 \to k[F] \to \bigoplus_{\facets\ F' \subset F} k[F'] \to \cdots \to
  \bigoplus_{\rays\ v \in F} k[v] \to k \to 0.
\end{equation}
The $\ZZ^d$-graded degree~$a$ part of this complex is zero unless $a \in
F$, in which case we get the relative chain complex $C_\spot(F,F')$,
where $a$~is in the relative interior of~$F'$.  The homology of such a
relative complex is zero unless $F' = F$.  Therefore, the only homology
of (\ref{eq:F}) is the canonical module $\omega_{k[F]}$, being the kernel
of the first map.%
\end{proof}

\begin{prop} \label{p:tot}
The total complex $\tot D(\Delta)$ of the Zeeman double complex is an
irreducible resolution of $k[Q]/I_\Delta$.
\end{prop}
\begin{proof}
The spectral sequence obtained by first taking vertical cohomology of
$D(\Delta)$ has $H_\partial D(\Delta) = E^1 = E^\infty$ by
Lemma~\ref{l:vert}.  The same lemma implies that the cohomology of
$\tot D(\Delta)$ is zero except in degree $p + q = 0$, and that the
nonzero cohomology has a filtration whose associated graded module is
$\bigoplus_{F \in \Delta} \omega_{k[F]}$.  On the other hand, the map
$k[Q] \to D(\Delta)$ sending $1 \mapsto \sum_{F \in \Delta} \epsilon_F F
\otimes e_F$ has kernel $I_\Delta$, for any choice of signs $\epsilon_F =
\pm 1$.  Choosing the signs
$$
  \epsilon_F = (-1)^{\dim F(\dim F + 1)/2} =
  \left\{
  \begin{array}{@{}r@{\ \ {\rm if}\ }l@{}}
	-1 & \dim F \equiv 1,2 \mod 4\\
	 1 & \dim F \equiv 0,3 \mod 4
  \end{array}
  \right.
$$
forces $(\delta + \partial)(\sum_{F \in \Delta} F \otimes e_F) = 0$,
thanks to the factor $(-1)^q$ in the definition~of~$\delta$.%
\end{proof}

\begin{cor} \label{c:irr}
Every summand in the minimal irreducible resolution of the quotient
$k[Q]/I$ by a reduced monomial ideal $I$ is isomorphic to $k[F]$ for some
face $F$ of~$Q$.
\end{cor}
\begin{proof}
Every summand in the total Zeeman complex of Proposition~\ref{p:tot} has
the desired form.  Now apply Theorem~\ref{t:irr}.\ref{split}.
\end{proof}

The spectral sequence in the proof of Proposition~\ref{p:tot} always
converges rather early, at~$E^1$.  The other spectral sequence, however,
obtained by first taking the horizontal cohomology~$H_\delta$, may be
highly nontrivial.

\begin{defn} \label{d:sequence}
The \bem{$\ZZ^d$-graded Zeeman spectral sequence} for the polyhedral
complex $\Delta$ is the spectral sequence $\ZZ E^\spot_{pq}(\Delta)$ on
the double complex $D(\Delta)$ obtained by taking horizontal homology
first, so $\ZZ E^2_{pq}(\Delta) = H_\partial H_\delta D(\Delta)$.
The \bem{ordinary Zeeman spectral sequence} for $\Delta$ is the
$\ZZ^d$-graded degree~$\0$ piece $\mathit{ZE}^\spot_{pq}(\Delta) = \ZZ
E^\spot_{pq}(\Delta)_\0$.
\end{defn}

\section{Characterization of Cohen--Macaulay quotients}

\label{sec:cm}

This section contains a characterization of Cohen--Macaulayness in terms
of irreducible resolutions coming from the Zeeman double complex
$D(\Delta)$.

\begin{defn} \label{d:cm}
The polyhedral complex $\Delta$ is \bem{Cohen--Macaulay over $k$} if the
local cohomology over $k$ of $\Delta$ near every face $G \in \Delta$
satisfies $H^i_G(\Delta) = 0$ for $i < \dim \Delta$.
\end{defn}

\begin{thm} \label{t:cm}
Let $I = I_\Delta$ be a radical monomial ideal.  The following are
equivalent.
\begin{enumerate}
\item \label{cm/k}
$\Delta$ is Cohen--Macaulay over $k$.

\item \label{ordinary}
The only nonzero vector spaces $\mathit{ZE}^1_{pq}(\Delta)$ lie in
column~$p = \dim(\Delta)$.

\item \label{E1}
The complex $\ZZ E^1(\Delta)$ is a minimal linear irreducible resolution
of~$k[Q]/I$.

\item \label{lin}
$k[Q]/I$ has a linear irreducible resolution.

\item \label{cm}
$k[Q]/I$ is a Cohen--Macaulay ring.
\end{enumerate}
\end{thm}
\begin{proof}
\ref{cm/k} $\Leftrightarrow$ \ref{ordinary}: The $\ZZ^d$-degree $\0$ part
of Lemma~\ref{l:hor} says that $\mathit{ZE}^1$ has $\bigoplus_{\dim G =
q} H^p_G(\Delta)$ in cohomological degree~$-q$.  The equivalence is now
immediate from Definition~\ref{d:cm}.

\ref{cm/k} $\implies$ \ref{E1}: The $E^1$ term in question is the complex
$H_\delta D(\Delta)$, with the differential $\partial$ in
Lemma~\ref{l:hor}.  That lemma together with Definition~\ref{d:cm}
implies that the horizontal cohomology $H_\delta D(\Delta)$ has one
column ($p = \dim \Delta$), that must therefore be a resolution of
something having the same Hilbert series as $k[Q]/I_\Delta$ by
Proposition~\ref{p:tot}.  Since $H^n_F(\Delta) = kF$ for facets $F \in
\Delta$, it is enough to check that the diagonal embedding $k[Q]/I_\Delta
\into \bigoplus_{\facets\ F \in \Delta} k[F]$ is contained in the kernel
of the first map of~$(H_\delta D(\Delta), \partial)$.  Suppose $\dim
\Delta = n$.  If $\dim G = n-1$ for some face $G \in \Delta$, then
\begin{eqnarray*}
  H^n_G(\Delta) &=& ({\textstyle\bigoplus kF)/\<\textstyle \sum
  \varepsilon(G,F) F\>},
\end{eqnarray*}
both sums being over all facets $F \in \Delta$ containing~$G$.  Now
calculate
$$
  \partial\Bigl(\sum_{\dim F = n\,} F \otimes e_F\Bigr)\ =\!\sum_{\dim F
  = n} \sum_{F \supset G} F \otimes \varepsilon(G,F) e_G\ =\!\!\!\!
  \sum_{\dim G = n-1} \Bigl(\sum_{F \supset G} \varepsilon(G,F) F\Bigr)
  \otimes e_G\ =\ 0.
$$

\ref{E1} $\implies$ \ref{lin}: Trivial.

\ref{lin} $\implies$ \ref{cm}: If $M$ is any module having a linear
irreducible resolution $\WW^\spot$ in which each summand of $\WW^i$ is
reduced, then $M$ is Cohen--Macaulay.  This can be seen by induction on
$d - \dim(M)$ via the long exact sequence for local cohomology $H^i_\mm$,
where $\mm$ is the graded maximal ideal.  The induction requires the
modules $\WW^i$ to be Cohen--Macaulay themselves, which holds because
$Q$~is saturated.

\ref{cm} $\implies$ \ref{cm/k}: $k[Q]/I$ being Cohen--Macaulay implies
that $\eext^i_{k[Q]}(M,\omega_{k[Q]})$ is zero for $i \neq d - n$, where
$n = \dim \Delta$.  In particular, if $\Omega^\spot$ is the $Q$-graded
part of the minimal injective resolution of $\omega_{k[Q]}$, then the
$i^\th$ cohomology of $\hhom(k[Q]/I,\Omega^\spot)$ is zero unless $i = d
- n$.  The complex $\Omega^\spot$ is the linear irreducible resolution of
$\omega_{k[Q]}$ in which each quotient $k[F]$ for $F \in \RR_{\geq 0} Q$
appears precisely once; see~(\ref{eq:F}).
Since $\hom(k[Q]/I,k[F]) = k[F]$ if $F \in \Delta$ and zero otherwise,
$\hhom(k[Q]/I,\Omega^\spot)$ is
$$
  0 \to \bigoplus_\twoline{F \in \Delta}{\dim F = n} k[F] \to \cdots \to
  \bigoplus_\twoline{F \in \Delta}{\dim F = \ell} k[F] \to \cdots \to k
  \to 0.
$$
If $a \in Q$ is in the relative interior of $G \in \Delta$, then the
$\ZZ^d$-graded degree~$a$ component of this complex is the homological
shift of $C^\spot(\Delta_G)$ whose $i^\th$ cohomology is
$H^{d-i}_G(\Delta)$.%
\end{proof}

\begin{remark}
Note the interaction of Theorem~\ref{t:cm} with the characteristic of
$k$: the horizontal cohomology of the Zeeman double complex $D(\Delta)$
can depend on $\char(k)$, just as the other parts of the theorem can.%
\end{remark}

When the semigroup $Q$ is $\NN^d$, so that $k[Q]$ is just the polynomial
ring in $d$~variables $z_1,\ldots,z_d$ over~$k$, the polyhedral complex
$\Delta$ becomes a simplicial complex.  Thinking of $\Delta$ as an order
ideal in the lattice $2^{[d]}$ of subsets of $[d] := \{1,\ldots,d\}$, the
\bem{Alexander dual} simplicial complex $\dv$ is the complement
of $\Delta$ in $2^{[d]}$, but with the partial order reversed.  Another
way to say this is that $\dv = \{[d] \minus F \mid F \not\in
\Delta\}$.

Theorem~\ref{t:cm} can be thought of as the extension to arbitrary normal
semigroup rings of the Eagon--Reiner theorem \cite{ER}, which concerns
the case $Q = \NN^d$, via the Alexander duality functors defined in
\cite{Mil2,Rom}.  To see how, recall that a $\ZZ$-graded
$k[\NN^d]$-module is said to have \bem{linear free resolution} if its
minimal $\ZZ$-graded free resolution over $k[z_1, \ldots, z_d]$ can be
written using matrices filled with linear forms.

\begin{cor}[Eagon--Reiner] \label{c:ER}
If $\Delta$ is a simplicial complex on $\{1, \ldots, d\}$, then $\Delta$
is Cohen--Macaulay if and only if $I_{\dv}$ has linear free
resolution.
\end{cor}
\begin{proof}
The minimal free resolution of $I_{\dv}$ is the functorial Alexander dual
(see \cite[Definition~1.9]{Rom} or \cite[Theorem~2.6]{Mil2} with
${\mathbf a} = \1$) of the minimal irreducible resolution of
$k[\NN^d]/I_\Delta$ guaranteed by Theorem~\ref{t:cm}.  Linearity of the
irreducible resolution translates directly into linearity of the free
resolution of~$I_\dv$.
\end{proof}

\begin{remark}
Is there a generalization of Theorem~\ref{t:cm} to the sequential
Cohen--Macaulay case that works for arbitrary saturated semigroups,
analogous and Alexander dual to the generalization \cite{HRW} of
Corollary~\ref{c:ER}?  Probably; and if so, it will likely say that the
ordinary and $\ZZ^d$-graded Zeeman spectral sequences collapse at $E^2$
(i.e.\ all differentials in $E^{\geq 3}$ vanish).
\end{remark}


\begin{remark}
The Alexander dual of the complex $\ZZ E^1(\Delta) = (H_\delta D(\Delta),
\partial)$, which provides a linear free resolution of $I_\dv$ in the
Cohen--Macaulay case, also provides the ``linear part'' of the free
resolution of $I_\dv$ when $\Delta$ is arbitrary \cite{RW}.  It is
possible to give an apropos proof of this fact using the Alexander dual
of the Zeeman spectral sequence for a Stanley--Reisner ring along with an
argument due to J.~Eagon \cite{Eag} concerning how to make spectral
sequences into minimal free resolutions.%
\end{remark}

\section{Remarks and further directions}

\label{sec:further}

Zeeman's original spectral sequence appears verbatim as the ordinary
Zeeman spectral sequence $\mathit{ZE}^\spot_{pq}$ in
Definition~\ref{d:sequence}, with $Q = \NN^d$.  Zeeman used his double
complex and spectral sequence to provide an extension of Poincar\'e
duality for singular triangulated topological spaces
\cite{zeeI,zeeII,zeeIII}.  When the topological space is a manifold, of
course, usual Poincar\'e duality results.  In the present context,
Zeeman's version of the Poincar\'e duality isomorphism should glue two
complexes of irreducible quotients of $k[\NN^d]$ together to form the
minimal irreducible resolution for the Stanley--Reisner ring of any
Buchsbaum simplicial complex---these simplicial complexes behave much
like manifolds.  This gluing procedure should work also for the more
general Buchsbaum polyhedral complexes~$\Delta$ obtained by considering
arbitrary saturated affine semigroups~$Q$.

Theorem~\ref{t:cm} is likely capable of providing a combinatorial
construction of the ``canonical \v Cech complex'' for $I_\dv$
\cite{Mil2,YanMonSup} when $\Delta$ is Cohen--Macaulay, or even Buchsbaum
(if the previous paragraph works).
Although $\dv$ has only been defined a~priori for simplicial complexes
$\Delta$, when $Q = \NN^d$, the definition of functorial squarefree
Alexander duality extends easily to the case of arbitrary saturated
semigroups.  The catch is that $I_\dv$ is not an ideal in~$k[Q]$, but
rather an ideal in the semigroup ring $k[Q^\star]$ for the cone $Q^\star$
dual to~$Q$.  Combinatorially speaking, the face poset of~$Q$ is not
usually self-dual, as it is when $Q = \NN^d$, so the process of
``reversing the partial order'' geometrically forces the switch to
$Q^\star$.  The functorial part of Alexander duality follows the same
pattern as the case $Q = \NN^d$: quotients $k[F]$ of $k[Q]$ are dual to
prime ideals $P_{F^\star}$ inside $k[Q^\star]$, where $F^\star$ is the
face of $Q^\star$ dual to~$F$.  This kind of construction is evident in
the work of Yanagawa \cite[Section~6]{YanMonSup}.

In general, irreducible resolutions---and perhaps other resolutions
by structure sheaves of subschemes---can be useful for computing the
$K$-homology classes of reduced subschemes that are unions of
transverally intersecting components.  When the ambient scheme is
regular, this method is an alternative to calculating free resolutions,
which produce $K$-\emph{co}homology classes.  In particular, this holds
for subspace arrangements in projective spaces.  This philosophy
underlies the application of irreducible resolutions in
\cite[Appendix~A.3]{grobGeom} to the definition of ``multidegrees''.

Note that when $Q \not\cong \NN^d$, irreducible resolutions are the only
finite resolutions to be had: free and injective resolutions of finitely
generated modules rarely terminate.  In particular, an understanding of
the Hilbert series of irreducible quotients of $k[Q]$---a polyhedral
problem---would give rise to formulae for Hilbert series of $Q$-graded
modules.  Similarly, algorithmic computations with irreducible
resolutions can allow explicit computation of injective resolutions,
local cohomology, and perhaps other homological invariants in the
$\ZZ^d$-graded setting over semigroup rings.


\footnotesize
\def\cprime{$'$}
\providecommand{\bysame}{\leavevmode\hbox to3em{\hrulefill}\thinspace}


\vbox{\footnotesize \baselineskip 10pt
\bigskip\noindent
%
%
Applied Mathematics, 2-363A, Massachusetts Institute of Technology, 77
Massachusetts Avenue,

Cambridge, MA 02139 \qquad email: {\tt ezra@math.mit.edu}}
\end{document}